\documentclass[10pt,a4paper]{article}

\usepackage[margin=1in]{geometry}
\usepackage{fourier}
\usepackage{amsmath,amsthm,amssymb}

\newtheorem{definition}{Definition}
\newtheorem{theorem}{Theorem}
\newtheorem{corollary}{Corollary}

\date{}

\begin{document}

\title{Reducing the Bias in Blocked Particle Filtering for High-Dimensional Systems \thanks{This work was supported by AFOSR/AOARD via AOARD-144042.}}

\author{Francesco Bertoli\thanks{F. Bertoli is with the Australian National University and NICTA. He is supported by NICTA.} \and Adrian N. Bishop\thanks{A.N. Bishop is with the Australian National University and NICTA. He is supported by NICTA and the Australian Research Council (ARC) via a Discovery Early Career Researcher Award (DE-120102873).}}

\maketitle

\begin{abstract}
Particle filtering is a powerful approximation method that applies to state estimation in nonlinear and non-Gaussian dynamical state-space models. Unfortunately, the approximation error depends exponentially on the system dimension. This means that an incredibly large number of particles may be needed to appropriately control the error in very large scale filtering problems. The computational burden required is often prohibitive in practice. Rebeschini and Van Handel (2013) analyse a new approach for particle filtering in large-scale dynamic random fields. Through a suitable localisation operation they reduce the dependence of the error to the size of local sets, each of which may be considerably smaller than the dimension of the original system. The drawback is that this localisation operation introduces a bias. In this work, we propose a modified version of Rebeschini and Van Handel's blocked particle filter. We introduce a new degree of freedom allowing us to reduce the bias. We do this by enlarging the space during the update phase and thus reducing the amount of dependent information thrown away due to localisation. By designing an appropriate tradeoff between the various tuning parameters it is possible to reduce the total error bound via allowing a temporary enlargement of the update operator without really increasing the overall computational burden. 
\end{abstract}

\section{Introduction}

Recursive Bayesian estimation (or filtering) is a technique for recursively estimating the state of a random process observed via noisy measurements.  If the underlying dynamical model is linear and Gaussian we have the celebrated Kalman filter \cite{kalmanfilter} which is an exact solution to the Bayesian filtering problem. Unfortunately, in many practical scenarios of interest, the Bayes filter is not exactly computable. Therefore, we seek techniques to approximate this ideal filter. The Kalman filter can be applied in more general settings \cite{kalmanextended,kalmanuscented} as an approximation. Particle filtering is a more general approximation method that is easily applied to nonlinear and non-Gaussian state-space models. The particle filter approximates the Bayesian filter via Monte Carlo simulation/sampling. The samples (or particles) are propagated through a sequential importance sampling mechanism that attempts to capture the dynamics of the unobservable process and the likelihood of the observations available. Other approximations exist such as Gaussian mixture filters etc. \cite{kalmanextended,anderson:79}.

The particle filter has been widely studied in theory and in countless practical applications \cite{gordon1993novel, kitagawa1996monte, ad-ndf-ng:01a}. In \cite{del2001stability} the authors prove that the error can be controlled uniformly in time, thus providing a solid mathematical support for application of the filter in numerous fields. Unfortunately, the particle filter computation is strongly dependent on the dimension of the underlying estimation problem. Specifically, the error bound grows exponentially with the system's dimension, making the filter infeasible in most high-dimensional applications. This problem is known as the \textit{curse of dimensionality} \cite{quang2010insight}. A heuristic explanation of this phenomenon for a particular case can be found in \cite{snyder2008obstacles}. In \cite{bickel2008sharp,quang2010insight} the authors give a precise relation between the dimension of the system and the number of particle required to avoid weight degeneracy \cite{ad-ndf-ng:01a}. The fact that the approximation error is exponential in the dimension and only inversely controlled by the sample size implies that an incredibly large number of particles are required when dealing with a high dimensional system if we want to control the error at a reasonable level. Obviously, a large number of particles means a heavy computational burden, that is often simply prohibitive. 

Recent studies \cite{beskos2011stability,beskos2011error,rebeschini2013can,rebeschini2013comparison} however suggest that high-dimensional particle filtering may be feasible in particular applications and/or if one is willing to accept a degree of systematic bias. In \cite{beskos2011stability}, the particle filter is applied in a static setting where the objective is to sample from some high-dimensional target distribution. In this case, through a sequence of intermediate and simpler distributions, it is shown that the particle filter will converge to a sampled representation of the target distribution with a typical Monte Carlo error (inverse in the number of particles) given a complexity on the order of the dimension squared. Although \cite{beskos2011stability} deals only, in essence, with a static problem of sampling from a fixed target distribution, the analysis introduces a novel way of thinking about high-dimensional particle filtering which may carry over to dynamic filtering problems. Related work appears in \cite{beskos2011error}.

\subsection{Background: The Motivating Paper}

In \cite{rebeschini2013can} the authors consider particle filtering in large-scale dynamic random fields. They assume the dynamics of the underlying process are localised to a neighbourhood of the field and the observations are local to each site. They exploit this idea by localising the algorithm during the update phase. They argue that the difficulty in high dimensional particle filtering is due largely to the dimension of the observation and the nonlinearity of the update operation. Therefore, they partition the field into independent blocks and correct every marginalised block separately. The posterior is simply the product of the blocked marginals. The real contribution of \cite{rebeschini2013can} is a descriptive and technical analysis that shows the error introduced due to the localisation procedure can be readily controlled if the dynamics of the random field at each site are only locally dependent on those sites within close proximity. The standard sampling approximation error is shown to be exponential in only the size of the individual blocks. The number of samples/particles controls the sampling approximation error at the typical rate while the error due to the localisation process is a systematic bias that can only be controlled through an increase in the block size. Since each block is updated independently, parallel implementation is readily applicable and the computational burden may be alleviated, albeit this remains to be seen in practice. While the results of \cite{rebeschini2013can} are at the proof-of-concept stage, the idea is incredibly powerful.

The authors in \cite{rebeschini2013can} show that although the total approximation error can be controlled uniformly in time, it suffers from a spatial inhomogeneity. Specifically, the nodes close to the block boundaries display a larger error than those far removed from the boundaries (as one might expect). A simple approach to average this spatial inhomogeneity is given in \cite{bertolibishop} where adaptive partitioning of the field is employed.

\subsection{Contribution}

In this paper we consider again the idea proposed in \cite{rebeschini2013can} and propose a modified particle filtering algorithm that displays an additional degree of freedom. The idea proposed herein is to enlarge the blocks during the update phase, allowing for more observations to be employed during the correction at each block. The main contribution is the addition of a new parameter that captures how much we enlarge each block prior to the update. Obviously, by enlarging each block prior to updating we reduce the bias error but we increase the complexity involved in updating each (enlarged) block. By designing an appropriate tradeoff between the various tuning parameters it is possible to reduce the total error bound via allowing a temporary enlargement of the update operator without increasing the overall computational burden.

\section{Problem Setup and Applications of the Blocked Filter}

We borrow the problem setup and notation directly from \cite{rebeschini2013can}. 

Consider a Markov chain $(X_n)_{n\geq 0}$ defined on a Polish state space $\mathbb{X}$ with transition density $p :\mathbb{X} \times \mathbb{X}\rightarrow \mathbb{R}$ with respect to a reference measure $\psi$. Moreover consider a process $(Y_n)_{n\geq 0}$, defined on a Polish space $\mathbb{Y}$, conditionally independent given $(X_n)_{n\geq 0}$, with a transition density $g:\mathbb{X}\times \mathbb{Y}\rightarrow \mathbb{R}$ with respect to a measure $\varphi$. The process $(X_n)_{n\geq 0}$ is observed via the process $(Y_n)_{n\geq 0}$.
Our aim is to estimate the probability of the state $X_n$ given the measurements up to that time and the initial condition $\mu$. Therefore we introduce the filter
$$
	\pi^{\mu}_n:= \mathbf{P}_\mu[X_n\in \cdot~ | Y_1,\cdots, Y_n]
$$
It can be easily seen, using Bayes rule, that the filter can be written in a recursive way
$$
	\pi_0^{\mu}=\mu,\qquad \pi^{\mu}_n=\mathsf{F}_n\pi^{\mu}_{n-1}
$$
where the operator $\mathsf{F}_n$ is defined as follows
$$
	(\mathsf{F}_n\rho)(A) :=\frac{\int \mathbf{1}_A(x) p(x_0,x) g(x,Y_n) \psi(dx)\rho(dx_0)}{\int  p(x_0,x) g(x,Y_n) \psi(dx)\rho(dx_0)}
$$
Moreover, the above operator is typically split into two sub-steps $\mathsf{F}_n=\mathsf{C}_n\mathsf{P}$ where
$$
	(\mathsf{P} \rho)(A) := \int \mathbf{1}_A(x) p(x_0,x)\psi(dx)\rho(dx_0)
$$
is a prediction step, and
$$
	(\mathsf{C}_n \rho)(A):= \frac{\int \mathbf{1}_A(x) g(x, Y_n) \rho(dx)}{\int g(x, Y_n) \rho(dx)}
$$
is a correction (or update) step. In the prediction step, the measure is transformed according to the density $p(\cdot,\cdot)$, while in the update step we use the new information $Y_n$ to correct the predicted measure. We then write the recursion as follows
$$
	\pi_{n-1}^\mu ~ \xrightarrow{\mathrm{prediction}} ~ \pi_{n-}^\mu := \mathsf{P}\pi_{n-1}^\mu ~ \xrightarrow{\mathrm{correction}} ~ \pi_n^\mu = \mathsf{C}_n\pi_{n-}^\mu
$$
The classic bootstrap particle filter uses $N$ particles (or samples) to approximate the measure $\pi_{n}^{\mu}$. Given a sampled approximation of $\pi_{n-1}^{\mu}$, the particles are first moved according to the transition $p(\cdot,\cdot)$ in order to approximate a sampled representation of the prediction. The update then computes a weighted posterior empirical measure via $g(\cdot,\cdot)$. Eventually, a resample step is added in order to avoid weight degeneracy \cite{ad-ndf-ng:01a}. More formally, denoting the bootstrap filter by $\hat{\pi}^{\mu}_n$, we have $\hat{\pi}^{\mu}_n=\hat{\mathsf{F}}_n\hat{\pi}^{\mu}_{n-1}$ where $\hat{\mathsf{F}}_n=\mathsf{C}_n\mathsf{S}^N \mathsf{P}$ and $\mathsf{S}^N$ represents the sampling operator here defined
$$
	\mathsf{S}^N\rho:=\frac{1}{N}\sum_{i=1}^N \delta_{x(i)}, \qquad x(i)~\mathrm{is~i.i.d.}\sim\rho
$$
It is possible to prove that
$$
	\sup_{|f|\leq 1} \mathbf{E}[\pi_n^{\mu}(f)-\hat{\pi}_n^{\mu}(f)]\leq {a_0}/{\sqrt{N}}
$$
with $a_0$ independent of time. Unfortunately, the constant $c$ typically depends (exponentially) on the dimension of the underlying problem. Intuition for this exponential dependence is given in \cite{rebeschini2013can, snyder2008obstacles}. 

We now consider the pair $(X_n,Y_n)$ as a random field $(X^v_n, Y^v_n)_{v\in V}$ indexed on a finite undirected graph $G=(V,E)$. The vertex set $V$ will represents the collection of sites and the edge set $E$ the spatial relationships between them. The cardinality of $V$ captures, in some sense, the dimension of interest. More formally, the spaces $\mathbb{X}$ and $\mathbb{Y}$ are defined as products $\mathbb{X}:=\prod_{v\in V} \mathbb{X}^v$, $\mathbb{Y}:=\prod_{v\in V} \mathbb{Y}^v$. The reference measures are products $\psi:=\bigotimes_{v\in V} \psi^v, \varphi:=\bigotimes_{v\in V} \varphi^v$, where $\psi^v$ and $\varphi^v$ are reference measures on $\mathbb{X}^v$ and $\mathbb{Y}^v$ respectively. The transition densities are defined as
$$
	p(x,z):=\prod_{v\in V} p^v(x,z^v), \qquad g(x,y):=\prod_{v\in V} g^v(x^v,y^v)
$$
where $p^v :\mathbb{X} \times \mathbb{X}^v\rightarrow \mathbb{R}$ and $g^v:\mathbb{X}^v\times \mathbb{Y}^v\rightarrow \mathbb{R}$ are densities with respect to the reference measures $\psi^v$ and $\varphi^v$. From the definition we can see that the observations $Y_n$ are assumed to be completely local, in the sense that $Y_n^v$ depends uniquely on the value assumed by $X_n^v$. The process $(X_n)_{n\geq 0}$ is local in the sense that the state at a site $v$ depends only on the state at nearby sites. We state this formally. Consider the graph $G$ equipped with the distance $d(v,v')$ defined by the number of hops along the shortest path connecting $v$ and $v'$. We can define the neighbourhood of a site $v$ as
$$
	N(v):=\{ v'\in V~ :~ d(v,v')\leq r\}
$$
where $r$ represents the range of interaction. Then we assume
$$
	p^v(x_1,z^v)=p^v(x_2,z^v) \quad \mathrm{whenever} \quad x_1^{N(v)}=x_2^{N(v)}
$$
where we write for $I\subseteq V$, $x^I=(x^i)_{i\in I}$. In other words, the random field $(X_n)_{n\geq 0}$ is local in the sense that given $X_0,\ldots, X_{n-1}$ the present state $X^v_n$ depends only on $X_{n-1}^{N(v)}$.

\subsection{Blocked Particle Filter}

In \cite{rebeschini2013can} the authors propose an application of the blocked filter algorithm to the field model just explained, exploiting the local dynamic dependencies. We briefly illustrate this algorithm. Consider a partition $\mathcal{K}=\{K_i\}_i$ of $V$ into non-overlapping blocks with a union equal to $V$. The idea is to create independence across blocks on $V$ by marginalising after the prediction step. We then update each block separately and finally we form $\pi^{\mu}_n$ via the product of the independent (updated) blocked marginals. More formally, consider the block operator $\mathsf{B}$ on the space $\mathcal{M}(\mathbb{X})$ of measures on $\mathbb{X}$, defined by
$$
	\mathsf{B}\rho := \bigotimes_{K\in \mathcal{K}} \mathsf{B}^K \rho
$$
where $\mathsf{B}^K\rho$ is the marginal of the measure $\rho$ on the subset $K\subseteq V$. Then the proposed block filter can be written as a recursion $\hat{\pi}^{\mu}_0=\mu$, $\hat{\pi}^{\mu}_n=\hat{\mathsf{F}}_n \hat{\pi}^{\mu}_{n-1}$ where the operator $\hat{F}_n=\mathsf{C}_n\mathsf{B}\mathsf{S}^N\mathsf{P}$ consists of four steps
$$
	\hat\pi_{n-1}^\mu ~ \xrightarrow{\mathrm{prediction/sampling}} ~ \hat\pi_{n-}^\mu = \mathsf{S}^N\mathsf{P}\hat\pi_{n-1}^\mu ~ \xrightarrow{\mathrm{blocking/correction}} ~ \hat\pi_n^\mu = \mathsf{C}_n\mathsf{B}\hat\pi_{n-}^\mu
$$

We make the following definition.

\begin{definition}
Given $\mu,\nu \in \mathcal{M}(\mathbb{X})$ and a subset $I\subseteq V$ we define a distance of the marginals on $I$ as follows
$$
	\VERT \mu - \nu \VERT_I := \sup_{f\in M(\mathbb{X}^I):|f|\leq 1} \mathbf{E}[|\mu(f)-\nu(f)|^2]^{\frac{1}{2}}
$$
where the expectation is taken with respect to the random sampling and $M(\mathbb{X}^I)$ is the class of measurable function on $\mathbb{X}$ that depends only on the values on $I$, that is $f(x)=f(y)$ when $x^I=y^I$. If $I=V$ we omit the subscript and write $\VERT \mu-\nu \VERT$. 
\end{definition}

With no expectation it follows that $\VERT\cdot\VERT$ is equivalent to the total variation which we write as $\|\cdot\|$. The two norms are interchangeable when no sampling occurs. 

Now, given a set $I\subseteq V$ we define the boundary and the interior 
$$
	\partial I := \{ v\in I ~|~ N(v) \not\subseteq I\},\qquad \mathrm{int}(I):=I\backslash\partial I
$$
and given a partition $\mathcal{K}$, we define the following quantities 
\begin{align*}
    	\Delta &:= \max_{v\in V} |N(v)| \\
	|\mathcal{K}|_\infty &:= \max_{K\in\mathcal{K}} |K| \\
	\Delta_\mathcal{K} &:= \max_{K\in\mathcal{K}} |\{K'\in\mathcal{K}:d(K,K')\le r\}|
\end{align*}
where the first quantity is independent of the partition. The result proven in \cite{rebeschini2013can} is the following.

\begin{theorem}[Blocked Particle Filter \cite{rebeschini2013can}]
\label{rebeschini main}
There exists a constant $0<\varepsilon_0<1$, depending only on the quantities $\Delta, \Delta_{\mathcal{K}}$ such that if there exists $\varepsilon_0<\varepsilon <1$ and $0<\kappa<1$ such that
$$
	\varepsilon< p^v(x, z^v) <\varepsilon^{-1}, \qquad \kappa<g^v(x^v,y^v)<\kappa^{-1} \qquad \forall x,z\in \mathbb{X},~ y \in \mathbb{Y},~ v\in V
$$
then for every $x\in \mathbb{X}$, $n\geq 0$, $K\in \mathcal{K}$ and $I\subseteq K$ we have
$$
	\VERT \pi^{\mu}_n-\hat{{\pi}}^{\mu}_n \VERT_I ~\leq~ \alpha \left[e^{-\beta_1d(I,\partial K)}+ \frac{e^{\beta_2 |\mathcal{K}|_{\infty}}}{\sqrt{N}}\right]
$$
where the constants $\alpha,\beta_1,\beta_2$ are positive, finite and dependent only on $\Delta, \Delta_{\mathcal{K}}, \varepsilon, \kappa$ and $r$.
\end{theorem}

The intuition is that the algorithm approximation error is exponential in $|K|$ rather then in $|V|$ but that the error at some individual locations increases with the proximity of those locations to the border of the blocks. This leads to a spatial inhomogeneity as seen in the first term of the bound. 

\subsection{Adaptively Blocked Particle Filter}

A first attempt to achieve a spatially homogeneous error bound can be found in \cite{bertolibishop}. The idea is to consider a finite number $m$ of partitions $\mathcal{K}_i$ and to apply them cyclically. Clearly we have to choose the partitions is such a way there is no node that is consistently close to a border. This condition is expressed by a bound on the average, or exponential average, of the border distance. Given $\beta>0$ write
$$
	\theta \leq \theta_m(v) = \frac{1}{m}\sum_{j=0}^{m-1} d(v,\partial K_j(v)) \quad \quad 0 < \phi_m(v) = \frac{1}{m}\sum_{j=0}^{m-1} e^{-\beta d(v,\partial K_j(v))} \leq \phi
$$

Clearly $\theta$ and $\phi$ represent how well balanced the collection of partitions are. Define $\Delta_d(v):=\max_s d(v,\partial K_s)$ and $\nabla_d(v):=\min_sd(v,\partial K_s)$.

\begin{theorem}[\cite{bertolibishop}]
\label{bertoli main}
There exists a constant $0<\varepsilon_0<1$, depending only on the quantities $\Delta$, $\Delta_{\mathcal{K}}$ such that if there exists $\varepsilon_0<\varepsilon <1$ and $0<\kappa<1$ such that
$$
	\varepsilon< p^v(x, z^v) <\varepsilon^{-1}, \qquad \kappa<g^v(x^v,y^v)<\kappa^{-1} \qquad \forall x,z\in \mathbb{X},~ y \in \mathbb{Y},~ v\in V
$$
then for every $x\in \mathbb{X}$, $n\geq 0$ and $v\in V$ we have
\begin{align*}
	\frac{1}{m}\sum_{k=0}^{m-1} \VERT \pi_{n-k}^\mu-\hat{\pi}_{n-k}^\mu \VERT_v &~\leq~ \alpha\bigg(\phi_m(v) + \frac{|\mathcal{K}|_{\infty}e^{\beta |\mathcal{K}|_{\infty}}}{\sqrt{N}}\bigg) \\
    & ~\leq ~ \alpha\bigg(e^{-\beta e^{-\beta(\Delta_d(v)-\nabla_d(v))} \frac{1}{m} \theta_m(v)}+\frac{|\mathcal{K}|_{\infty}e^{\beta |\mathcal{K}|_{\infty}}}{\sqrt{N}}\bigg)
\end{align*}
where $0<\alpha,\beta<\infty$ depend only on $\varepsilon$, $\kappa$, $r$, $\Delta$ and $|\mathcal{K}|_\infty := \max_s\max_{K\in\mathcal{K}_s} |K|$ in this case.
\end{theorem}

If $\theta = \theta_m(v) = \frac{1}{m}\sum_{j=0}^{m-1} d(v,\partial K_j(v))$ where $K_j(v)\in\mathcal{K}_j$ for all $v\in\mathcal{V}$, then the bound  is completely spatially invariant. See \cite{bertolibishop} for further discussion on this method.

\section{Enlarged Blocked Particle Filtering}

Suppose now we are given a partition $\mathcal{K}$ over $V$ but it turns out we are interested only in estimating the marginal of $\pi^{\mu}_n$ on a particular block $K\in\mathcal{K}$. We could first redefine the partition with a larger block encompassing $K$ and a bunch of single site blocks (to speed up the overall computation). It is of course not possible to define a partition in this manner for multiple blocks of interest. However, the idea proposed here is based on extending the state space by creating multiple independent copies of the measurements (and states) that are then used in different (and independent) enlarged blocks.

We introduce some new notation. Consider a parameter $b\geq 0$, that we will consider fixed throughout the rest of the paper. Then define, for any $K\in\mathcal{K}$, an enlarged block
$$
	\overline{K}:=\{v\in V ~|~ d(v,K)\leq b\}
$$
Now define the enlarged spaces
$$
	\mathbb{X}^E:=\prod_{K\in \mathcal{K}} \prod_{v\in \overline{K}} \mathbb{X}^v, \qquad \mathbb{Y}^E:=\prod_{K\in \mathcal{K}} \prod_{v\in \overline{K}} \mathbb{Y}^v
$$
Consider the collection $\overline{\mathcal{K}}=\{\overline{K}~:~ K\in \mathcal{K}\}$. This is no longer a partition of $V$. However, $\overline{\mathcal{K}}$ is a partition on $\mathbb{X}^E$, and here we can apply the blocking and updating operators associated with $\overline{\mathcal{K}}$. We use the superscript $E$ to note enlarged objects. The measures $\psi^E$ and $\varphi^E$ are defined straightforwardly. The block operator becomes
$$
	\mathsf{B}^E : \mathcal{M}(\mathbb{X})\rightarrow \mathcal{M}(\mathbb{X}^E), \quad\quad \mathsf{B}^E(\rho):=\bigotimes_{K\in \mathcal{K}} \mathsf{B}^{\overline{K}}\rho
$$
To update, we need the same operator $C_n$ redefined on the new space $\mathcal{M}(\mathbb{X}^E)$,
$$
	(\mathsf{C}^E_n\rho)(A) := \frac{\int\mathbf{1}_A(x)~\prod_{v\in \mathbb{X}^E}g^v(x^v,Y^v_n)~\rho(dx)}{\int \prod_{v\in \mathbb{X}^E}g^v(x^v,Y^v_n)~\rho(dx)}
$$
We also define
$$
	\mathsf{B}^{-1}:\mathcal{M}(\mathbb{X}^E)\rightarrow \mathcal{M}(\mathbb{X}), \quad\quad \mathsf{B}^{-1}(\rho):=\bigotimes_{K\in \mathcal{K}} \mathsf{B}^K\rho
$$
Now we can write the enlarged blocked filter algorithm as a recursion
$$
	\hat\pi_0^\mu=\mu, \qquad\qquad \hat\pi_n^\mu = \mathsf{\hat F}^E_n\hat\pi_{n-1}^\mu\quad(n\geq 1)
$$
where $\mathsf{\hat F}^E_n:=\mathsf{B}^{-1}\mathsf{C}^E_n\mathsf{B}^E\mathsf{S}^N\mathsf{P}$. Now we have five steps. Skipping the prediction/sample steps, graphically we have
$$
	\hat\pi_{n-}^\mu = \mathsf{S}^N\mathsf{P}\hat\pi_{n-1}^\mu ~ \xrightarrow{\mathrm{enlarging/blocking}} ~ \mathsf{B}^E\hat\pi_{n-}^\mu  ~ \xrightarrow{\mathrm{updating}} ~ \mathsf{C}^E_n\mathsf{B}^E\hat\pi_{n-}^\mu ~ \xrightarrow{\mathrm{marginalizing}} ~ \hat\pi_n^\mu
$$

To write out  the explicit expression of the filter we note that
$$
	(\mathsf{B}^{-1}\mathsf{C}^E_s\mathsf{B}^E\mathsf{P}\nu)(A)=(\mathsf{C}^E_s\mathsf{B}^E\mathsf{P}\nu)(A^E)
$$
where $A^E:=A\times (\mathbb{X}^E \backslash \mathbb{X}) $. Therefore, splitting a variable $z\in \mathbb{X}^E$ in $z=(x,z_E)$ with $x\in \mathbb{X}$ and $z_E\in \mathbb{X}^E\backslash\mathbb{X}$ (where now we put $E$ as subscript just for notational simplicity) and an enlarged block $\overline{K}=(K,K^E)$ where $K^E=\{v\in\overline{K}~:~ v\notin K\}$, we can write
$$
	(\mathsf{\tilde F}^E_s\nu)(A) = \frac{\int\mathbf{1}_{A^E}(z) \prod_{K'\in\mathcal{K}} \left[\prod_{w\in \overline{K'}}~p^w(x_0,z^w)~g^w(z^w,Y_s^w)~\nu(dx_0)\psi^{\overline{K'}}(dz^{\overline{K'}})\right]}
	{\int\prod_{K'\in\mathcal{K}}\left[\prod_{w\in \overline{K'}}~p^w(x_0,x^w)~g^w(x^w,Y_s^w)~\nu(dx_0)\psi^{\overline{K'}}(dz^{\overline{K'}}\right]}
$$
$$
	=\frac{\int\mathbf{1}_A(x)\prod_{K'\in\mathcal{K}} \left[\prod_{w\in K'}p^w(x_0,x^w)g^w(x^w,Y_s^w)\prod_{w\in K'^E}p^w(x_0,z_E^w)g^w(z_E^w,Y_s^w)\nu(dx_0)\psi^{K'}(dx^{K'})\psi^{K'^E}(dz_E^{K'^E})\right]}
	{\int\prod_{K'\in\mathcal{K}} \left[\prod_{w\in K'}p^w(x_0,x^w)~g^w(x^w,Y_s^w)\prod_{w\in K'^E}p^w(x_0,z_E^w)~g^w(z_E^w,Y_s^w)~\nu(dx_0)\psi^{K'}(dx^{K'})\psi^{K'^E}(dz^{K'^E}_E)\right]}
$$
where for $J\subseteq V$ we write $\psi^J(dx^J)=\prod_{v\in J} \psi^v(dx^v)$.

\section{Main Results and Discussion}

Define an ideal enlarged blocked filter $\tilde{\pi}^{\mu}_n=\tilde{\mathsf{F}}_n\dots\tilde{\mathsf{F}}_1\mu$ where $\tilde{\mathsf{F}}_s:=\mathsf{B}^{-1}\mathsf{C}^E_s\mathsf{B}^E\mathsf{P}$. Fix $I\subseteq V$. We then use the triangle inequality to decompose the error according to
$$
	\VERT \pi_n^\mu-\hat\pi_n^\mu \VERT_I ~\leq~ \VERT \pi_n^\mu-\tilde\pi_n^\mu \VERT_I + \VERT \tilde\pi_n^\mu-\hat\pi_n^\mu\VERT_I
$$
where we refer to the first and second decomposed terms as the bias and variance respectively. The bias represents the error introduced solely as a result of the blocking operation. In the standard bootstrap filter, this bias term vanishes and the typical analysis considers only the variance term. 

Going forward, we consider bounding both the bias and the variance. We stress however, that the bias is fundamentally more interesting as it pertains directly to the localisation idea considered herein. Indeed, the sampling operation that leads to the variance term could be replaced with other approximation techniques with no loss of generality (albeit a different approximation error than detailed subsequently). 

For sake of completeness/clarity we firstly state a result that includes both a bias and a variance bound.
\begin{theorem}[Main result]
\label{main result}
Suppose there exists a constant $0<\varepsilon_0<1$, depending only on $\Delta$ and $\Delta_{\mathcal{K}}$ and assume
$$
\varepsilon< p^v(x,z^v)<\varepsilon^{-1}, \quad \kappa<g^v(x^v,y^v)<\kappa^{-1} \qquad \forall x,z\in\mathbb{X},~ y\in\mathbb{Y},~ v\in V
$$
Then for every time $n\geq 0$, $x\in \mathbb{X}$, $K\in \mathcal{K}$ and $I\subseteq K$ we have
$$
	\VERT \pi^\mu_n - \hat{\pi}^\mu_n \VERT_I ~\leq~ \alpha\left[e^{-\beta_1d(I,\partial \overline{K})}+\frac{e^{\beta_2|\overline{\mathcal{K}}|_{\infty}}}{\sqrt{N}}\right]
$$
where the constants $0<\alpha,\beta_1,\beta_2<\infty$ depend only on $\varepsilon$, $\kappa$, $r$, $\Delta$, $\Delta_{\mathcal{K}}$, $\Delta_{\overline{\mathcal{K}}}$.
\end{theorem}

This single (total error) bound is derived in practice as two separate bounds which we now explicitly state.

\begin{theorem}[Bounding the Bias]
\label{bias theorem}
Assume there exists $0<\varepsilon<1$ such that
$$
	\varepsilon < p^v(x,z^v) < \varepsilon^{-1} \qquad \mathrm{for~all~} v\in V,~x,z\in\mathbb{X}
$$
and such that
$$
	\varepsilon >  \varepsilon_0 = \left(1-\frac{1}{18\Delta^2}\right)^{1/2\Delta}
$$
Let $\beta = -(2r)^{-1}\log 18\Delta^2(1-\varepsilon^{2\Delta})>0$. Then for every $n\geq 0$ we have
$$
	\VERT \pi_n^\mu-\tilde\pi_n^\mu \VERT_I ~\leq~ \frac{8e^{-\beta}}{1-e^{-\beta}}(1-\varepsilon^{2\Delta})|I|e^{-\beta d(I,\partial\overline{K})}
$$
for every $x\in\mathbb{X}$, $K\in\mathcal{K}$ and $I\subseteq K$.
\end{theorem}

The only difference between this bias bound and the bias bound in \cite{rebeschini2013can} is the presence of $|\overline{\mathcal{K}}|_{\infty}$ in place of $|\mathcal{K}|_{\infty}$. For a given partition $\mathcal{K}$ any enlargement of the blocks in $\mathcal{K}$ yielding $\overline{\mathcal{K}}$ results in a tighter bias bound as expected.

\begin{theorem}[Bounding the Variance]
\label{variance theorem}
Assume there exists $0<\kappa,\varepsilon<1$ such that
$$
	\varepsilon < p^v(x,z^v) < \varepsilon^{-1},\qquad \kappa<g^v(x^v,y^v)<\kappa^{-1} \qquad \mathrm{for~all~} v\in V,~x,z\in\mathbb{X},~y\in \mathbb{Y}
$$
and such that
$$
	\varepsilon >  \varepsilon_0 = \left(1-\frac{1}{6\Delta_{\overline{\mathcal{K}}}\Delta^2}\right)^{1/2\Delta}
$$
Let $\beta = -\log 6\Delta_{\overline{\mathcal{K}}}\Delta^2(1-\varepsilon^{2\Delta})>0$ where $\Delta_{\overline{\mathcal{K}}}:=\max|\{K\in\mathcal{K} : d(K,\overline{K})\leq r\}|$. Then for every $n\geq 0$ we have
$$
	\VERT \tilde\pi_n^\mu-\hat\pi_n^\mu\VERT_I ~\leq~ |I| \frac{64\Delta_{\mathcal{K}}}{1-e^{-\beta}} \frac{\varepsilon^{-4|\overline{\mathcal{K}}|_{\infty}}\kappa^{-4|\overline{\mathcal{K}}|_{\infty}\Delta_{\mathcal{K}}}}{\sqrt{N}}
$$
for every $x\in\mathbb{X}$, $K\in\mathcal{K}$ and $I\subseteq K$.
\end{theorem}

Again, the only significant difference between this variance bound and the variance bound in \cite{rebeschini2013can} is the presence of $|\overline{\mathcal{K}}|_{\infty}$ in place of $|\mathcal{K}|_{\infty}$. The variance depends inversely on the number of samples and exponentially in the size of the enlarged blocks.

\subsection{How to Use the Enlarged Blocked Filter}

Roughly, we now explain how one may implement the enlarged blocked filter to reduce the bias as compared with the algorithm proposed in \cite{rebeschini2013can} while maintaining a comparable variance and computational complexity.

Suppose firstly that one has a random field over $|V|$ sites and the computational power available (defining a bound on $N$) ensures that blocks of size $|V|/k$ can be readily handled for some $k>0$. Then the complexity of the blocked particle filter proposed in \cite{rebeschini2013can} can, in a sense, be regarded as being of order $\mathcal{O}(kN)$. Really, one can imagine $k$ particle filters running in parallel over each block and each with complexity on the order of $\mathcal{O}(N)$. 

To exploit the enlarged blocked particle filter, one should start with a larger number $c>k$ of smaller blocks which when enlarged are mostly of the size $|V|/k$. Then, the complexity of the enlarged blocked particle filter proposed herein is on the order $\mathcal{O}(cN)$. One immediately sees that the variance of the enlarged blocked particle filter is mostly on the same order as that of the algorithm proposed in \cite{rebeschini2013can} and the computational complexity has only increased linearly. However, in almost all cases (and certainly with well-designed partitions) one will achieve a reduction in the bias at any given site in the random field.

\subsection{Spatial Homogeneity}

We consider a special but interesting case in which a spatial homogeneous total error bound is obtained, the bias bound is better (tighter) than in \cite{rebeschini2013can}, and the computational requirements largely unchanged when compared with the algorithm in \cite{rebeschini2013can}. 

\begin{corollary}
Assume the same hypothesis of Theorem \ref{bias theorem}. Consider the partition $\mathcal{K}=\{v\}_{v\in V}$ and suppose $b>r$. Then for every $n\ge 0$, $x\in\mathbb{X}$, and $v\in V$ we have
$$
	\VERT \pi_n^\mu - \tilde\pi_n^\mu \VERT_v ~\leq~ \frac{8e^{-\beta}}{1-e^{-\beta}} (1-\varepsilon^{2\Delta}) e^{-\beta (b-r)}
$$
\end{corollary}

This bound is spatially homogeneous and with $b>r$ it is strictly less than the bias bound introduced in \cite{rebeschini2013can}. Note that while the bias bound here is spatially homogeneous, the actual bias may still be inhomogeneous since this result is potentially based on over bounding. On the other hand, it is possible to apply the adaptive scheme proposed in \cite{bertolibishop} with the enlarged blocked filter and potentially achieve true spatial homogeneity.

\subsection{Discussion on the Enlarged Blocked Filter}

The idea of the enlarged blocked particle filter is essentially based on the principle that larger blocks lead to a reduction in the bias introduced due to blocking. 

So, why not just start with larger blocks?

\begin{itemize}
	\item Well, irrespective of the size of the blocks, if one applies the standard blocked particle filter of \cite{rebeschini2013can} then there will always exist sites on the border of a block.
	\item If we extend (or enlarge) the blocks as proposed herein, we (typically) reduce the bias at each site (and particularly those sites that were on the border of a block in the original partition).
	\item If we increase the number of samples $N$ with a fixed number of larger blocks (in the original partition) then while we can reduce the variance we have no effect on the bias for those sites on the border.
	\item If we start with small blocks in the original partition and then simultaneously enlarge the blocks along with the number of samples $N$ then it may be possible maintain a given variance (or even reduce the variance) as compared to a partition with larger original block sizes but with a guaranteed smaller bias at each site.
\end{itemize}

The high-level point is that it is computationally more desirable to run a few extra parallel implementations of the particle filter (corresponding to more (enlarged) blocks) and obtain a tighter bias bound than it is to run a few less parallel implementations of the particle filter for the same variance bound but a larger bias bound. This is only possible through enlargement of the blocks as described herein.

Finally, we comment on the matter of consistency (as defined in say \cite{julier1997non}) and observational double counting. Consider the partition $\mathcal{K}=\{v\}_{v\in V}$ and suppose $\overline{K}=V$ for each $K=v\in\mathcal{K}$. Practically, following the standard prediction step, the enlarged blocked filter is of the form $\mathsf{B}^{-1}\mathsf{C}^E_n\mathsf{B}^E\rho$ which is mathematically equivalent to $\mathsf{B}\mathsf{C}_n\rho$. The point of this illustration is to highlight that even in this case, involving the most extreme enlargement possible, we are not double counting information or effectively applying measurements twice, and the enlarged blocked particle filter is consistent as per \cite{julier1997non}.

\subsection{Proof Strategy}

In this section we provide a summary of the proof strategy. Clearly the main result in Theorem \ref{main result} is immediately implied by Theorems \ref{bias theorem} and \ref{variance theorem}. Much of the technical analysis required in the proof of Theorem \ref{main result} is similar to that originally detailed in \cite{rebeschini2013can}. 

In the case of the bias $\VERT\pi_n^\mu-\tilde\pi_n^\mu\VERT_J$, one first derives a local stability property for the filter $\pi_n^\mu$ which implies that the marginal over a local set $J\subseteq V$ of the initial state $\mu$ is forgotten exponentially fast. Such a property also implies that any approximation errors in, say, the initial state are also forgotten. It then follows that if one can bound the one-step approximation error $\VERT\mathsf{F}_n\tilde\pi_{n-1}^\mu-\mathsf{\tilde F}_n\tilde\pi_{n-1}^\mu\VERT_J$ at any time, then in conjunction with the local stability property one will obtain a time-uniform bound on the bias over a local region of the field.

In the case of the variance $\VERT\tilde\pi_n^\mu-\hat\pi_n^\mu\VERT_J$, a similar idea is used except one first establishes stability for the ideal enlarged blocked filter $\tilde\pi_n^\mu$. Then, one must bound the one-step approximation error $\VERT\mathsf{\tilde F}_n\hat\pi_{n-1}^\mu-\mathsf{\hat F}_n\hat\pi_{n-1}^\mu\VERT_J$ at any time. Putting the stability property and the bound on the one-step approximation together, one achieves the desired time-uniform bound on the variance of a block in the adaptively blocked filter.

We have obviously glossed over much of the intricacies involved in the proof in this summary. For example, in the case of the bias, the property introduced in \cite{rebeschini2013can} and referred to as the decay of correlations must be established to hold uniformly in time for the ideal block filter $\tilde\pi_n^\mu$. This property captures a notion of spatial stability where the state at some site in the random field is forgotten as one moves away from that site. Rebeschini et al. provide a novel measure of this decay that allows them to establish local stability of the filter $\pi_n^\mu$ and to establish a bound on the one-step approximation error $\VERT\mathsf{F}_n\tilde\pi_{n-1}^\mu-\mathsf{\tilde F}_n\tilde\pi_{n-1}^\mu\VERT_J$. Conceptually, a property like the decay of correlations is necessary to establish such results. 

Summarising, the steps needed to prove the bias bound are
\begin{enumerate}
\item Proving a (local) stability result for the ideal Bayesian filter;
\item Proving that a desired decay of correlation property holds uniformly in time for the measure $\tilde{\pi}^x_n$;
\item Controlling the one time-step error introduced by the new enlarged blocked filter;
\item Putting all these results together and finalising Theorem \ref{bias theorem}.
\end{enumerate}

The variance analysis follows much the same path with the prime difficulty being establishment of local stability for the ideal enlarged blocked filter. Summarising the steps involved in proving the variance bound,
\begin{enumerate}
\item Proving a local stability result for the ideal enlarged blocked filter;
\item Controlling the one time-step error due to the sampling in the enlarged blocked particle filter;
\item Putting these results together and finalising Theorem \ref{variance theorem}
\end{enumerate}

\noindent The proof details are omitted in this version of the work due to their similarity with those details presented in \cite{rebeschini2013can}, but are available upon request.

\section{Concluding Remarks}

We have presented a modified version of the blocked particle filter originally proposed in \cite{rebeschini2013can}. The main feature of our algorithm is that we add a new parameter that can be tuned to decrease the bias as compared to \cite{rebeschini2013can}. The high-level argument for this approach is that it is computationally more desirable to run a few extra parallel implementations of the particle filter (corresponding to more (enlarged) blocks) and obtain a tighter bias bound than it is to run a few less parallel implementations of the particle filter for the same variance bound but a larger bias bound. This gain in bias reduction, with the same variance, and only a linear increase in the computational complexity, is only possible through enlargement of the blocks as described herein.

Finally, we also point out that the same adaptive approach to changing partitions proposed in \cite{bertolibishop} could be applied in the case of the enlarged blocked filter and this is an additional method for spatial smoothing and may be of interest in those cases in which the underlying model is time-varying. 

\bibliographystyle{plain}
\bibliography{}

\end{document}